\documentclass {amsart}

\usepackage{amsmath,amssymb,amsfonts,enumerate,amsthm, amscd,}
\usepackage[T1]{fontenc}
\usepackage[page]{appendix}

\newcommand{\Z}{\mathbb{Z}}
\newcommand{\F}{\mathbb{F}}
\newcommand{\N}{\mathbb{N}}
\newcommand{\lo}{\longrightarrow}

\newcommand{\fm}{\mathfrak{m}}

\newtheorem{thm}{Theorem}[section]
\newtheorem{cor}[thm]{Corollary}
\newtheorem{lem}[thm]{Lemma}

\newtheorem{defn}[thm]{Definition}

\newtheorem{exam}[thm]{Example}
\newtheorem{rem}[thm]{Remark}

\def\proof{{\parindent0pt {\bf Proof.\ }}}
\def\proofp{{\parindent0pt {\bf Proof of Theorem 2.3.\ }}}

\theoremstyle{definition}

\theoremstyle{remark}

\theoremstyle{Definition and Notation}

\begin{document}
\bibliographystyle{amsplain}


\title[On $p-$Rings] {On $p-$Ring}

\author{Mohammed Kabbour}
\address{Mohammed Kabbour\\Department of Mathematics, Faculty of Science and Technology of Fez, Box 2202, University S.M. Ben Abdellah Fez, Morocco.
$$ E-mail\ address:\ mkabbour@gmail.com$$}


\maketitle

\begin{abstract}
 In this paper, we introduced the concept of a $p$-ideal for a given ring. We provide necessary and sufficient condition for $\dfrac{R[x]}{(f(x))}$ to be a $p$-ring, where $R$ is a finite $p$-ring. It is also shown that the amalgamation of rings, $A\bowtie^fJ$ is a $p$-ring if and only if so is $A$ and $J$ is a $p$-ideal. Finally, we establish the transfer of this notion to trivial ring extensions.
 \end{abstract}
\bigskip


 \begin{section} {Introduction}

 \par All rings considered below are commutative with identity element
 $\neq0;$ and all modules are unital. Following N.H. McCoy and D. Montgomery \cite{McM}, a ring $R$ is said to be a $p$-ring ($p$ is a prime integer) if $x^p=x$ and $ px=0,$ for each $x\in R.$ Thus a Boolean ring, as a ring in which every element is idempotent, is simply 2-ring ($p=2$). Recall that a ring is said to be reduced
 if its nilradical is zero.\\
 \par The following conditions on a ring $R$ are equivalent:
 \begin{enumerate}
       \item For each $a$ in $R,$ there is some $b\in R$ such that $a=a^2b.$
       \item $R$ is a reduced ring and every prime ideal is maximal (i.e $R$ is a reduced 0-dimensional ring).
       \item For any maximal ideal $\fm$ of $R,$ the localization $R_\fm$ at $\fm$  is a field.
 \end{enumerate}
A ring satisfying the conditions as above is called a von Neumann regular ring. See for instance \cite{Gl,H}.\\
\par Let $A$ be a ring, $E$ an $A$-module and let $R=A\propto E$ be the set of pair $(a,e)$ with pairwise addition and multiplication is giving by $(a,e)(b,f)=(ab,af+be),\  R$ is called the trivial ring extension of $A$ by $E$ (also called the idealization of $E$ over $A$).
Considerable work, part of it summarized in Glaz's book \cite{Gl} and Huckaba's book \cite{H}, has been concerned with trivial ring extensions.\\
\par Let $A$ and $B$ be a pair of rings, $J$ an ideal of $B$ and let $f:A\lo B$ be a ring homomorphism. The following sub-ring of $A\times B:$ $$A\bowtie ^fJ=\{(a,f(a)+j)\ ;\ a\in A,j\in J\}$$ is said to be
 amalgamation of
 $A$ with $B$ along $J$ with respect to $f.$ Motivations and some applications of this construction, introduced by M. D'Anna, C.A. Finocchiaro and M. Fontana, are well  discussed with more detail  in the recent paper \cite{DFF}.\\
\par The main purpose of this paper is to give  new and original families of examples of $p$-rings. Also we investigate the transfer of this notion to trivial ring extensions and amalgamation of rings.
\bigskip
\bigskip
 \end{section}
\bigskip
\bigskip

 \begin{section} {Main results}

\bigskip

We state formally the definition of a $p$-ideal for a given ring.
\bigskip
\begin{defn}
Let $R$ be a ring and let $p$ be a prime integer. An ideal $I$ of $A$ is called a $p$-ideal if for each $x\in I:$ $$x^p=x \ \mbox{ and } \ px=0.$$
\end{defn}
\bigskip

From this definition, we can deduce that a ring $R$ is a $p$-ring if and only if every principal ideal of $R$ is a $p$-ideal.\\
 In the next theorem, we give a necessary and sufficient condition for $\Z/n\Z$ to have a nonzero $p$-ideal.
\bigskip

\begin{thm}\label{th2.2}
Let $n$ be a nonnegative integer and let $p$ be a prime integer.\begin{itemize}
                                                            \item If $v_p(n)=1$ ($v_p(n)$ is the $p$-valuation of $n$) then $\Z/n\Z$ has a unique nonzero $p$-ideal.
                                                            \item Otherwise $(0)$ is the unique $p$-ideal of $\Z/n\Z.$
                                                          \end{itemize}
\end{thm}
\proof
We say that every ideal of $\Z/n\Z$ has the form $k\Z/n\Z,$ where $k\in\{0,...,n-1\}.$ Assume that $v_p(n)=0$ and let $k$ be an integer such that $1\leq k\leq n-1.$ Suppose that $k\Z/n\Z$ is $p$-ideal, then $pkx\in n\Z$ for each element $x$ in $\Z.$ Thus $n$ divides $pk$ and so $n$ divides $k,$ witch is absurd. We deduce that $(0)$ is the unique $p$-ideal  of $\Z/n\Z.$\\
\par We shall need to use the following property:\\
Let $R_1,...,R_n$ be rings then every ideal of $R_1\times \cdots\times R_n$ has the form $I_1\times \cdots\times I_n,$ where $I_k$ is an ideal of $R_k$ for each $k\in\{1,...,n\}.$\\
On the  other hand, it is easy to see that $I_1\times \cdots\times I_n$ is a $p$-ideal if and only if so is $I_k$ for all $k\in\{1,...,n\}.$ Now suppose that $v_p(n)=1$ and let $q$ be the integer such that $pq=n.$ We denote $\F_p=\Z/p\Z,$ the Galois field of order $p.$  From the assumption (since $p$ is relatively prime to $q$), we can write  $\Z/n\Z \simeq \F_p \times\Z/q\Z.$ From the previous part of the proof $\F_p \times\Z/q\Z$ has a unique nonzero $p$-ideal, which is $\F_p\times (0).$\\
\par Assume that $\alpha=v_p(n)\geq 2.$ There is some positive integer $q,$ relatively prime to $p,$ such that $n=p^\alpha q.$ Hence  $\Z/n\Z \simeq \Z/p^\alpha\Z\times\Z/q\Z.$ Let $I$ be a $p$-ideal  of $\Z/p^\alpha\Z$. Since $px=0$ for each $x\in I,$ there exists some integer $k\in\{0,...,p-1\}$ such that $I=kp^{\alpha -1}\Z/p^\alpha\Z.$ From the assumption $p$ divides $kx\left((kxp^{\alpha-1})^{p-1}-1\right)$ for every integer $x.$ But $p$ is relatively prime to $(kxp^{\alpha-1})^{p-1}-1,$ thus $p$ must divides $k,$ therefore $k=0.$ We conclude that $I=(0),$ and so $(0)$ is the unique $p$-ideal of $\Z/n\Z.$
\qed\\

\bigskip
For example let $R$ be the ring $\Z/60\Z.$ Then we have, as follow, the list of all $p$-ideals of $R,$ where $p$ ranges over the set of prime integers:
                               \begin{itemize}
                                 \item $20\Z/60\Z$ is the unique  nonzero 3-ideal of $R.$
                                 \item $12\Z/60\Z$ is the unique  nonzero 5-ideal of $R.$
                                 \item $R$ has not a nonzero 2-ideal (since $v_2(60)=2$).
                                 \item $(0)$ is a $p$-ideal for each prime integer $p.$
                               \end{itemize}
\bigskip

\begin{thm}\label{th2.3}
Let $p$ be a prime integer and let $f(x)\in\F_p[x]$ be a nonconstant polynomial over $\F_p.$ Then $\dfrac{\F_p[x]}{(f(x))}$ contains a nonzero $p$-ideal if and only if $f(x)$ has at last one simple zero in $\F_p.$
\end{thm}
\bigskip
We need the following lemmas before proving Theorem\ref{th2.3}

\bigskip
\begin{lem}\label{lem2.4}
Let $R$ be a $p$-ring, $R[x]$ the polynomial ring over $R$ in the indeterminate $x$ and let $f(x)$ be an element of $R[x].$ Then $\dfrac{R[x]}{(f(x))}$ is a $p$-ring if and only if $f(x)$ divides $x^p-x.$
\end{lem}
\proof
We suppose that $\dfrac{R[x]}{(f(x))}$ is a $p$-ring. Then $(x+(f(x)))^p=x+(f(x)),$ therefore $x^p-x\in(f(x)).$ Conversely, assume that $f(x)$ divides $x^p-x$ and let $0\neq g(x)\in R[x].$ By induction on $n=\deg g,$ the degree of the polynomial $g(x),$ we claim that $\left(g(x)\right)^p=g\left(x^p\right).$ Indeed, it is certainly true for $n=0.$ Assume that the statement is true for each $k\leq n$ and that $\deg g=n+1.$ We put $g(x)=a_{n+1}x^{n+1}+g_1(x),$ where $0\neq a_{n+1}\in R$ and $g_1(x)\in R[x]$ such that $\deg g_1\leq n.$ By the binomial theorem, $$(g(x))^p=\left(a_{n+1}x^{n+1}\right)^p+\left(g_1(x)\right)^p=a_{n+1}^px^{p(n+1)}+g_1\left(x^p\right).$$
Thus $(g(x))^p=g\left(x^p\right),$ as desired.\\
On the other hand, $x^p-x$ divides $x^{kp}-x^k$ for each positive integer $k.$ Hence $x^p-x$ divides  $(g(x))^p-g(x),$ and so $(g(x)+(f(x)))^p=g(x)+(f(x)).$ Finally, it is easy to see that $p(g(x)+(f(x)))=0,$ so we have the desired result.
\qed\\
\bigskip

\begin{lem}\label{lem2.5}
Let $f(x)$ be an irreducible polynomial over $\F_p$ and let $k$ be a nonnegative integer. Then the following statements are equivalent:
\begin{enumerate}
  \item [{\rm(1)}] The ring $\dfrac{\F_p[x]}{\left(f^k(x)\right)}$ contains a nonzero $p$-ideal.
  \item [{\rm(2)}] $\dfrac{\F_p[x]}{\left(f^k(x)\right)}$ is a $p$-ring.
  \item [{\rm(3)}] $\dfrac{\F_p[x]}{\left(f^k(x)\right)}$ is isomorphic (as a ring) to $\F_p.$
\end{enumerate}
In this case $k=1$ and $\deg f=1.$
\end{lem}
\proof
(1) $\Longrightarrow$ (2): Let $I$ be a nonzero $p$-ideal of $\dfrac{\F_p[x]}{\left(f^k(x)\right)}.$ There is some $j\in\{0,...,k-1\}$ such that $I=\dfrac{f^j(x)\F_p[x]}{\left(f^k(x)\right)}.$ We get that $f^k(x)$ divides $$f^j(x)g(x)\left((f^j(x)g(x))^{p-1}-1\right),$$ for all $g(x)$ in $\F_p[x].$ Hence $f^{k-j}(x)$ divides $f^{j(p-1)}(x)-1$ (since $\F_p[x]$ is an integral domain). It follows that $j=0,$ since $f(x)$ is relatively prime with $f^m(x)-1$ for every nonnegative integer $m.$ We conclude that $I=\dfrac{\F_p[x]}{\left(f^k(x)\right)},$ as desired.\\

(2) $\Longrightarrow$ (3): By using the above lemma, we get that $f^k(x)$ divides $x^p-x.$ We denote $\F_p=\{a_0,...,a_{p-1}\}.$ For each $i\in\{0,...,p-1\},\ a_i$ is a root of the polynomial $x^p-x.$ Therefore $x^p-x=\left(x-a_0\right)\cdots \left(x-a_{p-1}\right).$ We conclude that $k=1$ and $f(x)=x-a_i$ for some $i$ in $\{0,...,p-1\},$ and so
$$\frac{\F_p[x]}{\left(f^k(x)\right)}=\frac{\F_p[x]}{\left(x-a_i\right)}\simeq \F_p$$\\

(3) $\Longrightarrow$ (1): Clear
\qed\\

\bigskip

\proofp Suppose that $f(x)=(x-a)g(x),$ where $a$ is an element of $\F_p$ and $g(x)\in \F_p[x]$ such that $g(a)\neq0.$ Then $$\frac{\F_p[x]}{(f(x))}\simeq \frac{\F_p[x]}{(x-a)}\times\frac{\F_p[x]}{(g(x))}\simeq\F_p\times\frac{\F_p[x]}{(g(x))}.$$
But $\F_p\times\dfrac{\F_p[x]}{(g(x))}$ has a nonzero $p$-ideal which is $\F_p\times(0).$ The sufficient condition is now straightforward.
\par Conversely, suppose that $\dfrac{\F_p[x]}{(f(x))}$ contains a nonzero $p$-ideal. We may assume that $f(x)$ is monic polynomial. Let $f(x)=f_1^{k_1}(x)...f_n^{k_n}(x)$ be the irreducible factors decomposition of $f(x)\ (f_i(x)$ is a monic irreducible polynomial and $k_i\in \N^*,$ for each $i\in\{1,...,n\}$). By applying Chinese remainder theorem, we deduce that $$\frac{\F_p[x]}{(f(x))}\simeq \frac{\F_p[x]}{\left(f_1^{k_1}(x)\right)}\times ... \times \frac{\F_p[x]}{\left(f_n^{k_n}(x)\right)}.$$
On the other hand, the finite product $I_1\times...\times I_n$ of ideals is a $p$-ideal if and only if so is $I_k$ for each $k\in\{1,...,n\}.$ We deduce that there exists $i\in\{1,...,n\}$ such that  $\dfrac{\F_p[x]}{\left(f_i^{k_i}(x)\right)}$ has a nonzero $p$-ideal. By Lemma \ref{lem2.5}, $k_i=1$ and $\deg f_i=1.$ This completes the proof of Theorem \ref{th2.3}.
\qed\\
\bigskip
\par Our next theorem is due to N.H. McCoy, for instance see \cite[Theorem 1]{Mc} in the case where $p=2,$ and \cite[Theorem 8]{Mc} in the general case. It is shown that any finite $p$-ring is isomorphic to a direct sum of copies of $\F_p.$ For the convenience of reader, we include here a sketch of the proof.
\bigskip

\begin{thm}
Let $R$ be a finite ring. Then $R$ is a $p$-ring with $n$ maximal ideals if and only if $R\simeq (\F_p)^n.$
\end{thm}
\proof
$\Longleftarrow )$ Since every finite direct product $R_1\times\cdots\times R_n$ of rings is $p$-ring if and only if so is $R_k$ for each $k\in\{1,...,n\},$ then $R$ is a $p$-ring. On the other hand, every maximal ideal of $R_1\times\cdots\times R_n$ has the form  $R_1\times\cdots\times R_{k-1}\times\fm_k\times R_{k+1}\times\cdots\times R_n,$ where $\fm_k$ is a maximal ideal of $R_k,$ and $k\in\{1,...,n\}.$ We denote $$J_k=\F_p\times\cdots\times\F_p\times(0)\times\F_p\times\cdots\times\F_p,$$ $(0)$ in its k$^{th}$ place and $\F_p$ elsewhere, for each $k\in\{1,...,n\}.$ Then $\{J_1,...,J_n\}$ is the set of all maximal ideals of $(\F_p)^n.$ We conclude that $R$ is a $p$-ring with $n$ maximal ideals.\\
\par $\Longrightarrow )$ Since $a^p=a$ for each $a$ in $R,$ then $R$ is a von Neumann regular ring. Therefore every prime ideal of $R$ is maximal and $R$ is a reduced ring. It follows that $\displaystyle \bigcap_{1\leq i\leq n}\fm_i =(0),$ where $\{\fm_1,...,\fm_n\}$ is the set of all maximal ideals of $R.$ By using Chinese remainder theorem we deduce that:
 $$R=\frac{R}{\fm_1\cap ...\cap\fm_n} \simeq \frac{R}{\fm_1}\times\cdots\times \frac{R}{\fm_n}.$$
Now, we need only shows that every $\dfrac{R}{\fm_k}$ is isomorphic to $\F_p.$ Let $k\in\{1,...,n\},$ we denote $\dfrac{R}{\fm_k}=\{a_1,...,a_q\}.$ By the assumption every field $\dfrac{R}{\fm_k}$ is a $p$-ring, then each $a_i$ is a root of the polynomial $x^p-x$ and so $q\leq p.$ But $\dfrac{R}{\fm_k}$ is a finite field of characteristic $p,$ then there exists a nonnegative integer $\alpha$ such that $q=p^\alpha.$ Thus $q=p$ and $\dfrac{R}{\fm_k}\simeq \F_p,$ completing the proof of the theorem.
\qed\\
\bigskip

\begin{rem}
Let $R$ be a semi local $p$-ring with $n$ maximal ideals. Then \begin{enumerate}
                                                                                          \item [{\rm(1)}] $R$ is a finite $p$-ring and has $p^n$ elements.
                                                                                          \item [{\rm(2)}] $R$ has $2^n$ ideals which are all $p$-ideals.
                                                                                        \end{enumerate}
\end{rem}
\proof Under the notations of the above proof, it suffices to show that $\dfrac{R}{\fm_k}$ is a finite field for each $k\in\{1,...,n\}.$ Since every element of $\dfrac{R}{\fm_k}$ is a root of the polynomial $x^p-x\in\dfrac{R}{\fm_k}[x],$ we have the required property. \qed\\

\bigskip
Now, we give a characterization that $\dfrac{R[x]}{(f(x))}$ is a $p$-ring, in the case when $R$ is a finite $p$-ring.
\bigskip

\begin{thm}\label{thm2.8}
Let $R$ be a finite $p$-ring and let $\{\fm_1,...,\fm_n\}$ be the set of all maximal ideals of $R.$ For every polynomial $f(x)$ in $R[x]$ and $j\in\{1,...,n\},$ we denote by $f_j(x)$ the reduction of $f(x)$ modulo $\fm_j$ i.e $f_j(x)=\displaystyle\sum_{i=0}^k\left(a_i+\fm_j\right)x^i \in\frac{R}{\fm_j}[x],$ where $\displaystyle f(x)=\sum_{i=0}^ka_i x^i.$ Then $\dfrac{R[x]}{(f(x))}$ is a $p$-ring if and only if for each $j\in\{1,...,n\},\ f_j(x)$ splits with distinct roots in the field $\dfrac{R}{\fm_j}.$
\end{thm}
\proof
Under the above hypothesis, we get that $R[x]\simeq\dfrac{R}{\fm_1}[x]\times ...\times \dfrac{R}{\fm_n}[x],$  as a ring, via the map: $g(x)\mapsto (g_1(x),...,g_n(x)),$ where $g_j(x)$ is the reduction of $g(x)$ modulo $\fm_j.$ For each $j\in\{1,...,n\},$ we put $R_j=\dfrac{R}{\fm_j}.$ Then $R_j$ is a $p$-ring. On the other hand, the map $$\varphi:R[x]\lo \frac{R_1[x]}{(f_1(x))}\times ...\times \frac{R_n[x]}{(f_n(x))}$$ defined by $\varphi(g(x))=\left(g_1(x)+(f_1(x)),...,g_n(x)+(f_n(x))\right)$ is a surjective ring homomorphism. Also we have the following equality $\ker\varphi=(f(x)).$  Thus $\dfrac{R[x]}{(f(x))}$ is isomorphic to $\dfrac{R_1[x]}{(f_1(x))}\times ...\times \dfrac{R_n[x]}{(f_n(x))}.$ It follows that $\dfrac{R[x]}{(f(x))}$ is a $p$-ring if and only if so is $\dfrac{R_j[x]}{(f_j(x))},$ for each $j\in\{1,...,n\}.$ Now we can apply Lemma \ref{lem2.4} to prove that $\dfrac{R_j[x]}{(f_j(x))}$ is a $p$-ring if and only if $f_j(x)$ has $\deg f_j$ distinct roots in $R_j.$ This completes the proof of Theorem \ref{thm2.8}.
\qed\\

\bigskip
The next example illustrates the above results.
\bigskip

\begin{exam}\rm
Let $p$ be a prime integer of the form $8n+1,$ for some non negative integer $n.$ Consider the polynomial over $\F_p^4$ defined by $$f(x)=(1,-1,2,-2)+(0,0,1,1)x^2+(1,1,0,0)x^n.$$ Then $\dfrac{\F_p^4[x]}{(f(x))}$ is a finite $p$-ring with $2n+4$ maximal ideals and $p^{2n+4}$ elements.
\end{exam}

\proof
Under the above notations, we have
$$f_1(x)=1+x^n,\ f_2(x)=-1+x^n,\ f_3(x)=2+x^2 \mbox{ and } f_4(x)=-2+x^2.$$
 It is easy to see that $f_1(x)$ and $f_2(x)$ divide $x^{p-1}-1,$ hence $f_1(x)$ and $f_2(x)$ split with distinct roots in $\F_p.$ Also $x^{4n}+1$ divides $x^p-x,$ then $x^{4n}+1$ splits. Let $a$ be a root of the polynomial $x^{4n}+1\in\F_p[x],$ hence $\left(\dfrac{a^{2n}+1}{a^n}\right)^2=2$ and $\left(\dfrac{a^{2n}-1}{a^n}\right)^2=-2.$ Therefore $f_3(x)$ and $f_4(x)$ have distinct zeros  in $\F_p,$ since $f_3'(a)=f_4'(a)\neq 0\ , (f_j'(x)$ is the derivative of $f_j(x)$). The result then follows from Theorem \ref{thm2.8}.
\qed\\
\bigskip

In the next theorem we give our main result about the transfer of
$p$-ring property to amalgamation of rings.

\bigskip

\begin{thm}
Let $A$ and $B$ be a pair of rings, $J$ an ideal of $B,\ f:A\lo B$ a ring homomorphism and let $A\bowtie^fJ$ be the amalgamation of $A$ with $B$ along $J$ with respect to $f.$ Then $A\bowtie^fJ$ is a $p$-ring if and only if so is $A$ and $J$ is a $p$-ideal of $B.$
\end{thm}
\proof
$\Longrightarrow)$
Let $a\in A$ and $j\in J.$ It is easy to see that $(a,f(a))^p=(a^p,f(a^p))$ and $(0,j)^p=(0,j^p).$ But $(a,f(a))^p=(a,f(a))$ and  $(0,j)^p=(0,j),$ then $a^p=a$ et $j^p=j.$ Obviously $pa=0$ and $pj=0,$ since $p(a,f(a)+j)=0.$  We have the desired implication.\\

$\Longleftarrow)$ Assume that $A$ is a $p$-ring and $J$ is a $p$-ideal of $B.$ Let $(a,f(a)+j)\in A\bowtie^fJ.$ By the binomial theorem (which is valid in any commutative ring), $$(a,f(a)+j)^p=\left(a^p,f(a^p)+j^p+\sum_{k=1}^{p-1}
\binom{p}{k}
j^kf(a^{p-k})\right).$$
Since $j^kf(a^{p-k})\in J$ and $p$ divides $\dbinom{p}{k},$
for each $k\in\{1,...,p-1\},$ then $$(a,f(a)+j)^p=\left(a^p,f(a^p)+j^p\right)=(a,f(a)+j).$$ On the other hand, $p(a,f(a)+j)=(pa,f(pa)+pj)=0.$ It follows that $A\bowtie^fJ$ is a $p$-ring.
\qed\\

\bigskip
\begin{exam}\rm
Let $A$ be the set of all sequences of elements of $\F_p$ and let $B=\Z/n\Z,$ with $n=p(p+1).$ By using Theorem\ref{th2.2}, the principal ideal $(p+1)B$ is a $p$-ideal of $B.$ Consider the mapping $f:A\lo B$ defined by $f(a)=(p+1)a_0,$ where $a=\left(a_k+p\Z\right)_{k\in\N}.$ It is easy to see that $f$ is a ring homomorphism. On the other hand, the set all functions of a non empty set $X$ into a $p$-ring is also a $p$-ring. Hence $A$ is a $p$-ring. From the above theorem $A\bowtie^f(p+1)B$ is a $p$-ring.
\end{exam}
\bigskip

The following corollary is an immediate consequence of the above theorem.

\begin{cor}
Let $A$ be a ring and let $A\bowtie I$ be the amalgamated duplication  of $A$ along an ideal $I$ of $A.$ Then $A\bowtie I$ is a $p$-ring if and only if so is $A.$
\end{cor}
\bigskip

We end this paper by giving a necessary and sufficient condition for the trivial ring extension, $A\propto E,$ to be a von Neumann regular ring (resp., a $p$-ring).
\bigskip
\begin{thm}
Let $A$ be a ring, $E$ an $A$-module and let $A\propto E$ be the trivial ring extension of $A$ by $E.$ Then
 \begin{enumerate}
   \item [{\rm(1)}]  $A\propto E$ is a von Neumann regular ring if and only if so is $A$ and $E=\{0\}.$
   \item [{\rm(2)}]  $A\propto E$ is a $p$-ring if and only if so is $A$ and $E=\{0\}.$
 \end{enumerate}

\end{thm}
\proof
(1) We say that every maximal ideal of $A\propto E$ has the form $\fm\propto E,$ where $\fm$ is a maximal ideal of $A.$ Let $M$ be a maximal ideal of $A\propto E.$ By \cite[Theorem 4.1]{AW}, $(A\propto E)_M\simeq A_\fm\propto E_\fm,$ where $M=\fm\propto E.$ Thus $(A\propto E)_M$ is a field if and only if so is $A_\fm$  and $E_\fm=\{0\}.$ We deduce that $A\propto E$ is a von Neumann regular ring if and only if so is $A$ and $E_\fm=\{0\},$ for all maximal ideal $\fm$ of $A.$ We have the desired result.\\

\par (2) Assume that $A\propto E$ is a $p$-ring. It is easy to see that every sub-ring of $p$-ring is also a $p$-ring. It follows that $A$ is a $p$-ring. On the other hand, $E=\{0\}$ since $A\propto E$ is a von Neumann regular ring. We can also deduce this result from the following equalities: $$(a,x)=(a,x)^p=(a^p,pa^{p-1}x)=(a^p,0),$$ since $(a,x)^n=(a^n,na^{n-1}x)$ for every nonnegative integer $n,$ and $p(b,y)=0$ for every element $(b,y)$ of $A\propto E.$\\

The sufficient condition is obvious.
\qed\\


\bigskip
\noindent {\bf Acknowledgements.} The author thank the referee
for his/her careful reading of this work.

\bigskip

 \end{section}
\bigskip\bigskip

\end{document}